\documentclass[11pt, a4paper]{article}
\usepackage[english]{babel}
\usepackage[utf8]{inputenc}
\usepackage[T1]{fontenc}
\usepackage{a4wide}
\usepackage[pdftex]{graphicx}
\usepackage{verbatim}
\usepackage[width=0.8\textwidth]{caption}
\usepackage{amsmath,amssymb,amsopn,bm,dsfont}
\usepackage{float,placeins,flafter,longtable,array}
\usepackage{color}
\usepackage{natbib}
\usepackage{booktabs}
 \makeatletter
\newtheorem{thm}{Theorem}[section]

\newtheorem{lemma}{Lemma}[section]
\newtheorem{prop}[thm]{Proposition}
\newtheorem{definition}[thm]{Definition}

\newtheorem{hyp}{Assumption}
\numberwithin{equation}{section}

\newcommand{\sgn}{\text{sgn}}

\newcommand{\R}{\mathbb{R}}
\newcommand{\N}{\mathbb N}

\newcommand{\indep}{\perp \!\!\! \perp}
\newcommand{\cov}{\text{cov}}

\newcommand{\convNor}[1]{\stackrel{d}{\longrightarrow} \mathcal{N}\left(#1\right)}
\newcommand{\convD}{\stackrel{d}{\longrightarrow}}
\newcommand{\convP}{\stackrel{p}{\longrightarrow}}
\newcommand{\Vh}{\widehat{V}}

\renewcommand{\section}{\@startsection{section}{2}{0mm}{-1.5\baselineskip}{1\baselineskip}{\normalfont\large\bfseries}}
\renewcommand{\subsection}{\@startsection{subsection}{2}{0mm}{-1.2\baselineskip}{1\baselineskip}{\normalfont\normalsize\bfseries}}
\renewcommand{\subsubsection}{\@startsection{subsubsection}{3}{0mm}{-0.8\baselineskip}{0.4\baselineskip}{\normalfont\normalsize\itshape}}

\date{}

\begin{document}

\title{Empirical MSE Minimization to Estimate a Scalar Parameter}

\author{Cl\'{e}ment de Chaisemartin\thanks{University of California at Santa Barbara, clementdechaisemartin@ucsb.edu%
} \and Xavier D'Haultf\oe{}uille%
\thanks{de Chaisemartin: University of California at Santa Barbara (email: clementdechaisemartin@ucsb.edu); D'Haultf\oe uille: CREST-ENSAE (email: xavier.dhaultfoeuille@ensae.fr)}}

\maketitle ~\vspace{-1cm}

\begin{abstract}

We consider the estimation of a scalar parameter, when two estimators are available.  The first is always consistent. The second is inconsistent in general, but has  a smaller asymptotic variance than the first, and may be consistent if an assumption is satisfied. We propose to use the weighted sum of the two estimators with the lowest estimated mean-squared error (MSE). We show that this third estimator dominates the other two from a minimax-regret perspective: the maximum asymptotic-MSE-gain one may incur by using this estimator rather than one of the other estimators is larger than the maximum asymptotic-MSE-loss.
\end{abstract}



\section{Introduction}

We consider the estimation of a scalar parameter $\beta_0$ when two estimators are available. The first, is $\sqrt{n}-$consistent. The second is inconsistent in general, but it has a smaller asymptotic variance than the first, and it may be $\sqrt{n}-$consistent if the data generating process satisfies an assumption $H_0$. Hereafter, those two estimators are respectively referred to as the consistent and efficient estimators. 

\medskip
To fix ideas, we consider two of the many examples where this set-up we is applicable. In stratified randomized experiments, the parameter of interest is the average treatment effect (ATE). To estimate it, one may use the propensity score matching estimator \citep[see][]{hirano2003efficient}, which is $\sqrt{n}$-consistent and asymptotically normal under some assumptions. Alternatively, one may regress the outcome of interest on strata fixed effects and units' treatment status, and use the coefficient of the treatment in that regression. It follows from, e.g., Equation (3.3.7) in \cite{angrist2008mostly} that this estimator is $\sqrt{n}$-consistent and asymptotically normal for a weighted average of the effect of the treatment in each strata. One can also show that under some assumptions, the asymptotic variance of the strata fixed effects estimator is smaller than that of the propensity score matching estimator. So if the treatment effect is constant across strata, the strata fixed effects estimator is $\sqrt{n}$-consistent for the ATE, and it is more efficient than the propensity score matching estimator. But if the treatment effect is heterogeneous, the strata fixed effects estimator is inconsistent.

\medskip
Another example where our set-up is applicable is a linear and constant treatment effect model, where the treatment is potentially endogenous, but one has an instrumental variable at hand. Then, the 2SLS estimator is $\sqrt{n}-$consistent for the treatment effect. On the other hand, the OLS estimator is only $\sqrt{n}-$consistent if the treatment is actually exogenous, but its asymptotic variance is smaller than that of the 2SLS estimator.

\medskip
To estimate $\beta_0$, we propose to use $\widehat{\beta}_{MSE}$, the weighted sum of the consistent and efficient estimators with the lowest estimated mean-squared error (MSE). We show that this third estimator dominates the other two from a minimax-regret perspective: the maximum asymptotic-MSE-gain one may incur by using $\widehat{\beta}_{MSE}$ rather than one of the other estimators is larger than the maximum asymptotic-MSE-loss that one may incur by doing so.

\medskip
We also consider a family of pre-test estimators $(\widehat{\beta}_{MSE,\lambda})_{\lambda\geq 0}$, where $\lambda$ indexes the critical value used in the pre-test. First, we test whether the consistent and efficient estimators are equal. If the test is accepted, $\widehat{\beta}_{MSE,\lambda}$ is equal to the efficient estimator. If the test is rejected, $\widehat{\beta}_{MSE,\lambda}$ is equal to a convex combination of the consistent and efficient estimators. We show that such estimators have similar properties as $\widehat{\beta}_{MSE}$. However, $\widehat{\beta}_{MSE}$ dominates all of them from a minimax-regret perspective.

\medskip
We then extend the initial result by considering situations where one has two estimators at hand: one is $r_n-$consistent, where $r_n/\sqrt{n}\rightarrow 0$, and the other is inconsistent in general, but may be $\sqrt{n}-$consistent if the data generating process satisfies an assumption $H_0$. Such situations may for instance arise in regression discontinuity (RD) designs. Then, non-parametric estimators such as the one proposed by \cite{hahn2001identification} are $n^{2/5}-$consistent for the average treatment effect at the cut-off under weak conditions. On the other hand, the estimator using, say, linear regressions to the left and to the right of the cut-off without restricting the sample to observations in a narrow bandwidth around the cut-off is $\sqrt{n}-$consistent if the potential outcomes' CEFs are indeed linear in the running variable, but inconsistent otherwise. Again, we show that $\widehat{\beta}_{MSE}$ dominates the $r_n-$consistent estimator from a minimax-regret perspective, under mild assumptions.

\medskip
The idea of combining consistent and inconsistent, or unbiased and biased  estimators of a parameter has a long tradition in statistics and econometrics. \cite{green1991james} have proposed an estimator related to $\widehat{\beta}_{MSE}$, in the context of a normal model with known variances. Their estimator is related to the shrinkage estimator in \cite{charles1956inadmissibility} and \cite{james1961estimation}, and when the parameter of interest is of dimension greater than three, its MSE is lower than that of the unbiased estimator of the parameter of interest. \cite{judge2004semiparametric} and \cite{mittelhammer2005combining} have proposed an estimator similar to that in \cite{green1991james}, without assuming normality or that the variances are known, and study its asymptotic MSE, again when the parameter of interest is of dimension greater than three. \cite{cheng2019uniform} have considered an estimator related to $\widehat{\beta}_{MSE}$, in a GMM context with some valid and some potentially misspecified moment conditions. They show that asymptotically, the MSE of their estimator is uniformly smaller than that of the GMM estimator using only the valid moment conditions. Again, they focus on a multivariate parameter with a dimension greater than four. 
In the context of a linear model with at least three endogenous variables and instruments, \cite{hansen2017stein} proposes to use a weighted average of the OLS and 2SLS estimators, with weights that depend on the Hausman-Wu statistic in a test of equality between the OLS and 2SLS estimators. He shows that this estimator has a lower asymptotic MSE than that of the 2SLS estimator. Finally, \cite{breusch2011fixed} have considered $\widehat{\beta}_{MSE}$ in a panel data context, in a case where the parameter of interest is univariate. However, they do not study the theoretical properties of that estimator.

\medskip
To the best of our knowledge, our paper is the first to study the theoretical properties of $\widehat{\beta}_{MSE}$ when the parameter of interest is univariate. This case is particularly relevant for policy evaluation and treatment choice. Even when one measures the effect of the policy on several outcomes, one is ultimately interested in summarizing those effects into a monetary assessment of the benefits of the policy, to be compared to its cost. Contrary to the previous literature, we do not find that the MSE of $\widehat{\beta}_{MSE}$ dominates uniformly that of the consistent estimator. However, we show that $\widehat{\beta}_{MSE}$ dominates the other two estimators from a minimax-regret perspective, thus giving a theoretical justification to its use to estimate a univariate parameter. Our paper also seems to be the first to consider the combination of estimators with different rates of convergence, which may be relevant in a number of contexts, such as RD designs.

\medskip
The remainder of the paper is organized as follows. Section 2 presents the set up and main results. Section 3 presents some extensions. Section 4 presents the proofs of the results.

\section{Set up and main results}

We are interested in a parameter $\beta_0\in \R$. To estimate $\beta_0$, we use a sample of size $n$. $\widehat{\beta}_C$ is $\sqrt{n}-$consistent and asymptotically normal for $\beta_0$. $\widehat{\beta}_E$ is $\sqrt{n}-$consistent and asymptotically normal for $\beta_E$. In general, $\beta_E\ne \beta_0$, but under an assumption on the data generating process $H_0$, $\beta_E=\beta_0$. The asymptotic variance of $\widehat{\beta}_E$ is smaller than that of $\widehat{\beta}_C$, so under $H_0$, $\widehat{\beta}_E$ is a more efficient estimator of $\beta_0$ than $\widehat{\beta}_C$. Moreover, the asymptotic variance of $\widehat{\beta}_E$ is smaller than that of any weighted sum of $\widehat{\beta}_E$ and $\widehat{\beta}_C$, which implies that the asymptotic covariance of $\widehat{\beta}_E$ and $\widehat{\beta}_C$ is equal to the asymptotic variance of $\widehat{\beta}_E$. Finally, we have estimators $\widehat{V}(\widehat{\beta}_C)$, $\widehat{V}(\widehat{\beta}_E)$, and $\widehat{\cov}(\widehat{\beta}_C,\widehat{\beta}_E)$ of the variances of $\widehat{\beta}_C$ and $\widehat{\beta}_E$ and of their covariance, that are such that $n\widehat{V}(\widehat{\beta}_C)$, $n\widehat{V}(\widehat{\beta}_E)$, and $n\widehat{\cov}(\widehat{\beta}_C,\widehat{\beta}_E)$ are consistent for their asymptotic variances and covariances. We summarize these conditions in Assumption \ref{hyp:setup} below:
\begin{hyp}\label{hyp:setup} (Set-up)
\begin{enumerate}
\item \label{hyp:setup_p1_asnorm} We have $$\sqrt{n}\left( \begin{array}{c}
\widehat{\beta}_C-\beta_0 \\
\widehat{\beta}_E-\beta_E
\end{array}
\right)\convNor{0,\left( \begin{array}{cc}
\sigma^2_C & \sigma^2_E  \\
\sigma^2_E & \sigma^2_E
\end{array}\right)},$$
with $\sigma^2_C<\sigma^2_E$.
\item \label{hyp:setup_p2_consistent_var_estimation} $n\widehat{V}(\widehat{\beta}_C)\convP\sigma^2_C$, $n\widehat{V}(\widehat{\beta}_E)\convP\sigma^2_E$, and $n\widehat{\cov}(\widehat{\beta}_C,\widehat{\beta}_E)\convP \sigma^2_E$.
\end{enumerate}
\end{hyp}
We consider the following estimator of $\beta_0$:
\begin{definition}\label{definition_betaMSE} (The empirical-MSE-minimizing estimator of $\beta_0$)

\medskip
Let
\begin{equation}\label{eq:betaMSE}
\widehat{\beta}_{MSE}=\widehat{p}\widehat{\beta}_E+\left(1-\widehat{p}\right)\widehat{\beta}_C,
\end{equation}
where
\begin{equation}\label{eq:popt}
\widehat{p}=\arg\min_{p\in\R} ~~p^2\left(\widehat{\beta}_E-\widehat{\beta}_C\right)^2+p^2\widehat{V}(\widehat{\beta}_E)+(1-p)^2\widehat{V}(\widehat{\beta}_C)+2p(1-p)\widehat{\cov}(\widehat{\beta}_C,\widehat{\beta}_E).
\end{equation}
\end{definition}
$\widehat{\beta}_{MSE}$ is the weighted sum of $\widehat{\beta}_E$ and $\widehat{\beta}_C$ with the lowest estimated mean-squared error (MSE).
Solving the problem in Equation \eqref{eq:popt} yields
\begin{equation}\label{eq:popt_closedform}
\widehat{p}=\frac{\widehat{V}(\widehat{\beta}_C)-\widehat{\cov}(\widehat{\beta}_C,\widehat{\beta}_E)}{\left(\widehat{\beta}_E-\widehat{\beta}_C\right)^2+\widehat{V}(\widehat{\beta}_E-\widehat{\beta}_C)}.
\end{equation}

Theorem \ref{thm_main} below gives the asymptotic distribution of $\widehat{\beta}_{MSE}$ and compares it with that of $\widehat{\beta}_C$. In particular, we compare the MSE of the asymptotic distribution of the two estimators.\footnote{If the second moments of the normalized estimators converge, our results provide a comparison of the asymptotic MSE of the estimators.  To avoid the issue that convergence in distribution does not imply convergence in $L^2$, one could consider instead, as \cite{cheng2019uniform}, a winsorized version of the square loss.}

\begin{thm}\label{thm_main} (Asymptotic distribution of $\widehat{\beta}_{MSE}$)

Suppose Assumption \ref{hyp:setup} holds.
\begin{enumerate}
	\item If $\beta_0\ne\beta_E$ and $\beta_0$ and $\beta_E$ do not depend on $n$, 	$\sqrt{n}\left(\widehat{\beta}_{MSE}-\beta_0\right)\convNor{0,\sigma^2_C}$.
	\item If $\beta_0=\beta_E$, $\sqrt{n}\left(\widehat{\beta}_{MSE}-\beta_0\right)\convD U_0$, where $U_0$ is such that $E(U_0)=0$ and $V(U_0)\in (\sigma^2_E, \sigma^2_C)$.
	\item If $\beta_E=\beta_0+h/\sqrt{n}$ for some $h\in \R$, 	$\sqrt{n}\left(\widehat{\beta}_{MSE}-\beta_0\right)\convD U_h$, where $(U_h)_{h\in \R}$ is such that
	$$\max_{h\in \R}\left[E(U_h^2)-\sigma^2_C\right]<\max_{h\in \R}\left[\sigma^2_C-E(U_h^2)\right].$$
\end{enumerate}
\end{thm}

\medskip
Point 1 of Theorem \ref{thm_main} shows that if $\beta_0\ne\beta_E$, $\widehat{\beta}_{MSE}$ and $\widehat{\beta}_{C}$ have the same asymptotic distribution. On the other hand, Point 2 shows that if $\beta_0=\beta_E$, their asymptotic distributions differ, and the MSE of the asymptotic distribution of $\widehat{\beta}_{MSE}$ is larger than that of $\widehat{\beta}_E$ but smaller than that of $\widehat{\beta}_{C}$. This comes from the fact that under $H_0$, $\widehat{p}$ converges in distribution to a nongenerate distribution. In other words, $\widehat{p}$ does not perform consistent ``model selection'': it does not converge to 1 if $H_0$ holds and to 0 otherwise.

\medskip
The asymptotic approximations under fixed values of $\beta_0$ and $\beta_E$ may not give good approximations of the finite sample behavior of the estimators. Instead, in Point 3 of the theorem, we compare the MSE of their asymptotic distributions under alternatives local to $\beta_0=\beta_E$, namely $\beta_E=\beta_0+h/\sqrt{n}$. We find that under this type of asymptotics, the MSE of the asymptotic distribution of $\widehat{\beta}_{MSE}$ is not always smaller than that of $\widehat{\beta}_{C}$. This phenomenon is reminiscent of Hodges' estimator, whose asymptotic distribution has a smaller MSE than that of the standard estimator if the true parameter is 0,  but whose maximal risk increases without bound as $n\rightarrow\infty$ \citep[see, e.g.][pp. 440-443]{lehmann1998}. An important difference with Hodges' estimator is that here, the maximum asymptotic-MSE-gain one may incur by using $\widehat{\beta}_{MSE}$ rather than $\widehat{\beta}_{C}$ is larger than the maximum asymptotic-MSE-loss. Thus, $\widehat{\beta}_{MSE}$ dominates $\widehat{\beta}_{C}$ from a minimax regret perspective \citep[see][]{savage1951theory}. It is straightforward to show that $\widehat{\beta}_{MSE}$ also dominates $\widehat{\beta}_{E}$ from a minimax regret perspective.

\section{Extensions}

\subsection{A family of ``pre-test'' estimators} 
\label{sub:modified_versions_of_the_estimator}

When $\beta_0=\beta_E$, $\widehat{\beta}_{MSE}$ is not equivalent to $\widehat{\beta}_E$, and it has a higher asymptotic variance. This is because the estimated squared bias $(\widehat{\beta}_E-\widehat{\beta}_C)^2$ in \eqref{eq:popt} includes some noise and is not negligible even if $\beta_0=\beta_E$. We now consider a modified version of $\widehat{\beta}_{MSE}$ that uses a smaller estimator of the squared bias. Specifically, we replace \eqref{eq:popt} by
\begin{align*}
\widehat{p}_{\lambda}=\arg\min_{p\in\R} \quad & p^2\max\left[0,\left(\widehat{\beta}_E-\widehat{\beta}_C\right)^2 - \lambda \Vh(\widehat{\beta}_E-\widehat{\beta}_C) \right] \\
& + p^2\widehat{V}(\widehat{\beta}_E)+(1-p)^2\widehat{V}(\widehat{\beta}_C)+2p(1-p)\widehat{\cov}(\widehat{\beta}_C,\widehat{\beta}_E)
\end{align*}
for some $\lambda\geq 0$. We then let $\widehat{\beta}_{MSE,\lambda}=\widehat{p}_{\lambda} \widehat{\beta}_E + (1-\widehat{p}_{\lambda}) \widehat{\beta}_C$.
$\widehat{\beta}_{MSE,\lambda}$ can be viewed as a pre-test estimator. Assume that one uses $\widehat{V}(\widehat{\beta}_E)$ to estimate $\cov(\widehat{\beta}_C,\widehat{\beta}_E)$. Then, let $F_1$ denote the cdf of a $\chi^2_1$ distribution and let $\alpha=F_1(\lambda)$. To compute, $\widehat{\beta}_{MSE,\lambda}$, one first runs a level-$\alpha$ test of $H_0$, using the fact that $(\widehat{\beta}_E-\widehat{\beta}_C)^2/\Vh(\widehat{\beta}_E-\widehat{\beta}_C)\convD \chi^2_1$. If $H_0$ is accepted, then $\widehat{\beta}_{MSE,\lambda}=\widehat{\beta}_E$. If $H_0$ is rejected, $\widehat{\beta}_{MSE,\lambda}$ is equal to a convex combination between $\widehat{\beta}_E$ and $\widehat{\beta}_C$, where the weight assigned to $\widehat{\beta}_C$ depends on how far we are from accepting $H_0$.
\begin{prop}\label{prop:reduced_bias}
Suppose Assumption \ref{hyp:setup} holds and let $\lambda$ be a positive real number.
\begin{enumerate}
	\item If $\beta_0\ne\beta_E$ and $\beta_0$ and $\beta_E$ do not depend on $n$, $\sqrt{n}\left(\widehat{\beta}_{MSE,\lambda}-\beta_0\right)\convNor{0,\sigma^2_C}$.
	\item If $\beta_0=\beta_E$, $\sqrt{n}\left(\widehat{\beta}_{MSE,\lambda}-\beta_0\right)\convD U_{0,\lambda}$, where $U_{0,\lambda}$ is such that $E(U_{0,\lambda})=0$, $\lambda \mapsto V(U_{0,\lambda})$ is strictly decreasing and $\lim_{\lambda\to+\infty}  V(U_{0,\lambda})=\sigma^2_E$.
	\item If $\beta_E=\beta_0+h/\sqrt{n}$ for some $h\in \R$, $\sqrt{n}\left(\widehat{\beta}_{MSE,\lambda}-\beta_0\right)\convD U_{h,\lambda}$, where $(U_{h,\lambda})_{h\in \R}$ is such that for all $\lambda\neq 0$
\begin{equation}
\max_{h\in \R}\left[E(U_{h,\lambda}^2) - E(U_{h,0}^2) \right]>\max_{h\in \R}\left[E(U_{h,0}^2) - E(U_{h,\lambda}^2)\right],		
		\label{eq:ineq_minimax}
\end{equation}
\end{enumerate}
\end{prop}

Points 1 and 2 in Proposition \ref{prop:reduced_bias} are similar to those in Theorem \ref{thm_main}, with the additional point that under $H_0$, the asymptotic, quadratic risk of $\widehat{\beta}_{MSE,\lambda}$ decreases and gets closer to that of $\widehat{\beta}_E$ as $\lambda$ increases. However, the third point shows that from a minimax regret perspective, such estimators are dominated by our intial estimator $\widehat{\beta}_{MSE}$. The reason is that the decrease of $\lambda\mapsto E(U_{0,\lambda}^2)$ does not compensate for the quick increase of  $\lambda \mapsto \max_{h\in \R}E(U_{h,\lambda}^2) - E(U_{h,0}^2)$.
This echoes the discussion in \cite{leeb2005}: as we move closer to an estimator based on a consistent model selection, the maximal asymptotic risk increases without bound.


\subsection{Averaging estimators with different rates of convergence} 
\label{sub:averaging_estimators_with_different_rates_of_convergence}

In this subsection, we assume that $\widehat{\beta}_C$ is $r_n-$consistent for $\beta_0\in \R$, for some sequence $(r_n)_{n\in\N}$ such that $r_n\to\infty$, $r_n/n^{1/2}\to 0$. The estimator $\widehat{\beta}_E$ is still $\sqrt{n}-$consistent and asymptotically normal for $\beta_E$, which may be equal to $\beta_0$ under an assumption on the data generating process $H_0$. We also assume we have estimators $\widehat{V}(\widehat{\beta}_C)$, $\widehat{V}(\widehat{\beta}_E)$, and $\widehat{\cov}(\widehat{\beta}_C,\widehat{\beta}_E)$ of the variances of $\widehat{\beta}_C$ and $\widehat{\beta}_E$ and of their covariance, that are such that $r_n^2\widehat{V}(\widehat{\beta}_C)$, $n\widehat{V}(\widehat{\beta}_E)$, and $n^{1/2} r_n\widehat{\cov}(\widehat{\beta}_C,\widehat{\beta}_E)$ are consistent for their asymptotic variances and covariances. We summarize these conditions in Assumption \ref{hyp:setup_nonparam} below:
\begin{hyp}\label{hyp:setup_nonparam} (Set-up)
\begin{enumerate}
\item \label{hyp:setup_p1_asnorm_nonparam} There exists a sequence $(r_n)_{n\in\N}$, $r_n\to\infty$ and $r_n/n^{1/2}\to 0$, such that $$\left( \begin{array}{c}
r_n\left(\widehat{\beta}_C-\beta_0\right) \\
\sqrt{n}\left(\widehat{\beta}_E-\beta_E\right)
\end{array}
\right)\convNor{\left( \begin{array}{c}
\mu  \\
0
\end{array}
\right),\left( \begin{array}{cc}
\sigma^2_C & \rho  \\
\rho & \sigma^2_E
\end{array}
\right)}.$$
\item \label{hyp:setup_p2_consistent_var_estimation_nonparam} $r_n^2\widehat{V}(\widehat{\beta}_C)\convP\sigma^2_C$, $n\widehat{V}(\widehat{\beta}_E)\convP\sigma^2_E$, and $n^{1/2} r_n\widehat{\cov}(\widehat{\beta}_C,\widehat{\beta}_E)\convP \rho$.
\end{enumerate}
\end{hyp}

\medskip
Theorem \ref{thm_main_nonparam} below gives the asymptotic distribution of $\widehat{\beta}_{MSE}$ defined in Equation \ref{eq:betaMSE} above, under Assumption \ref{hyp:setup_nonparam} rather than Assumption \ref{hyp:setup}.

\begin{thm}\label{thm_main_nonparam} (Asymptotic distribution of $\widehat{\beta}_{MSE}$ under Assumption \ref{hyp:setup_nonparam})
Suppose Assumption \ref{hyp:setup_nonparam} holds.
\begin{enumerate}
\item If $\beta_0\ne\beta_E$ and $\beta_0$ and $\beta_E$ do not depend on $n$, $r_n\left(\widehat{\beta}_{MSE}-\beta_0\right)\convNor{0,\sigma^2_C}$.
\item If $\beta_0=\beta_E$, $r_n\left(\widehat{\beta}_{MSE}-\beta_0\right)\convD U_0$, where $U_0$ is such that $E(U_0^2)<E(V^2)$.
\item If $\beta_E=\beta_0+h/r_n$ for some $h\in \R$ and $|\mu/\sigma_C|\leq 0.4$, then $r_n\left(\widehat{\beta}_{MSE}-\beta_0\right)\convD U_h$, where $(U_h)_{h\in \R}$ is such that
$$\max_{h\in \R}\left[E(U_h^2)-(\mu^2+\sigma^2_C)\right]<\max_{h\in \R}\left[\mu^2+\sigma^2_C- E(U_h^2)\right].$$
    \end{enumerate}
\end{thm}

\medskip
Point 1 of Theorem \ref{thm_main_nonparam} shows that if $\beta_0\ne\beta_E$, $\widehat{\beta}_{MSE}$ and $\widehat{\beta}_{C}$ have the same asymptotic distribution. On the other hand, Point 2 shows that if $\beta_0=\beta_E$, their asymptotic distributions differ, and the MSE of the asymptotic distribution of $\widehat{\beta}_{MSE}$ is smaller than that of $\widehat{\beta}_{C}$. In Point 3 of the theorem, we compare the MSE of their asymptotic distributions under alternatives local to $\beta_0=\beta_E$, namely $\beta_E=\beta_0+h/r_n$. Under this type of asymptotics, the maximum asymptotic-MSE-gain one may incur by using $\widehat{\beta}_{MSE}$ rather than $\widehat{\beta}_{C}$ is larger than the maximum asymptotic-MSE-loss, provided the first-order bias is no greater in absolute value than 0.4$\sigma_C$. Again, $\widehat{\beta}_{MSE}$ dominates $\widehat{\beta}_{C}$ from a minimax regret perspective, provided the first-order bias of $\widehat{\beta}_{C}$ is not too large.


\section{Proofs} 
\label{sec:proofs}

We use the folowing lemma below.

\begin{lemma}\label{lem:expect_symm}
	Suppose that $f$ is an odd function such that $f(x)>0$ for all $x>0$ and $Z$ has an even density $g$ that is strictly decreasing on $\R^+$. Then, $\sgn(E[f(x+Z)])=\sgn(x)$ for all $x\in\R$.
\end{lemma}

\textbf{Proof:} the result holds if $x=0$, because $E[f(Z)]=E[f(-Z)]=-E[f(Z)]$. For any $x<0$,
$E[f(Z+x)]=E[f(-Z+x)]=E[-f(Z-x)]$, so it suffices to show that $E[f(Z+x)]>0$ for $x>0$. We have
\begin{align*}
	E[f(Z+x)] & =\int_{\R} f(z+x)g(z)dz\\
	& = \int_{\R^+} f(z)g(z-x)dz + \int_{\R^-} f(z)g(z-x)dz \\
	& =  \int_{\R^+} f(z)g(z-x)dz - \int_{\R^-} f(-z)g(x-z)dz \\
	& =  \int_{\R^+} f(z)[g(z-x) - g(z+x)]dz.
\end{align*}
Now, for all $z\in (0,x]$, $|z-x|=x-z < x+z$ so $g(z-x) > g(z+x)$. If $z>x$, $|z-x|=z-x < z+x$ so again, $g(z-x) > g(z+x)$. The result follows since $f(z)>0$ on $(0,\infty)$.

\subsection{Theorem \ref{thm_main}} 
\label{sub:theorem}

\subsubsection*{Proof of Point 1}

If $\beta_0\ne\beta_E$,  it follows from \eqref{eq:popt_closedform}, Assumption \ref{hyp:setup} and the continuous mapping theorem that
$$n\widehat{p}\convP \frac{\sigma^2_C- \sigma^2_E}{(\beta_E-\beta_0)^2}.$$
Moreover,
\begin{equation}\label{eq:diff_MSE_C}
\widehat{\beta}_{MSE}-\widehat{\beta}_C=\widehat{p}\left(\widehat{\beta}_E-\widehat{\beta}_C\right).
\end{equation}
Hence, by the continuous mapping theorem again, $$n\left(\widehat{\beta}_{MSE}-\widehat{\beta}_C\right)\convP\frac{\sigma^2_C- \sigma^2_E}{\beta_E-\beta_0}.$$
Therefore,
$$\sqrt{n}\left(\widehat{\beta}_{MSE}-\widehat{\beta}_C\right)=o_P(1).$$

\subsubsection*{Proof of Point 2}

Let $(V,W)$ be a normal vector with mean $(0,0)$, variances $(\sigma^2_C,\sigma^2_C-\sigma^2_E)$, and covariance $-\left(\sigma^2_C-\sigma^2_E\right)$. Let
\begin{equation}\label{eq:defU}
U_0=V+W\frac{\sigma^2_C-\sigma^2_E}{W^2+\sigma^2_C-\sigma^2_E}.
\end{equation}
\begin{align*}
&\sqrt{n}\left(\widehat{\beta}_{MSE}-\beta_0\right)\\
=&\sqrt{n}\left(\widehat{\beta}_{C}-\beta_0\right)\\
+&\frac{n\widehat{V}(\widehat{\beta}_C)-n\widehat{\cov}(\widehat{\beta}_C,\widehat{ \beta}_E)}{\left(\sqrt{n}\left(\widehat{\beta}_E-\beta_E-\left(\widehat{\beta}_C-\beta_0\right)\right)\right)^2+n\widehat{V}(\widehat{\beta}_E-\widehat{\beta}_C)}\sqrt{n}\left(\widehat{\beta}_E-\beta_E-\left(\widehat{\beta}_C-\beta_0\right)\right)\\
\convD&~U_0.
\end{align*}
The first equality follows from Equations \eqref{eq:diff_MSE_C} and  \eqref{eq:popt_closedform} and from $\beta_E=\beta_0$. The convergence in distribution arrow follows from Assumption \ref{hyp:setup}, the Slutsky lemma, and the continuous mapping theorem.

\medskip
$E(U_0)=0$ as $\phi: w\mapsto w\frac{\sigma^2_C-\sigma^2_E}{w^2+\sigma^2_C-\sigma^2_E}$ is such that $\phi(-w)=-\phi(w)$ and the pdf of $W$ is symmetric around 0.

\medskip
Let $\Psi=V+W$. The vector $(\Psi,W)$ is normally distributed with $E(\Psi)=0$ and $\cov(\Psi,W)=0$. Hence, $\Psi\indep W$ and we have
$$U_0 = \Psi + W\left(\frac{\sigma^2_C - \sigma^2_E}{W^2+\sigma^2_C-\sigma^2_E} -1 \right),$$
Hence, $V(U_0)>V(\Psi)=\sigma^2_E$. Moreover,
\begin{align*}
V(U_0)-\sigma^2_C=&E\left(U_0^2\right)-E\left(V^2\right)\\
=&E((U_0-V)(U_0+V))\\
=&E\left(W\frac{\sigma^2_C-\sigma^2_E}{W^2+\sigma^2_C-\sigma^2_E}\left(2V+W\frac{\sigma^2_C-\sigma^2_E}{W^2+\sigma^2_C-\sigma^2_E}\right)\right)\\
=&E\left(W\frac{\sigma^2_C-\sigma^2_E}{W^2+\sigma^2_C-\sigma^2_E}\left(-2W+2\Psi+W\frac{\sigma^2_C-\sigma^2_E}{W^2+\sigma^2_C-\sigma^2_E}\right)\right)\\
=&E\left(W^2\frac{\sigma^2_C-\sigma^2_E}{W^2+\sigma^2_C-\sigma^2_E}\left(-2+\frac{\sigma^2_C-\sigma^2_E}{W^2+\sigma^2_C-\sigma^2_E}\right)\right)\\
<&0.
\end{align*}
The first equality follows from $E(U_0)=E(V)=0$, the third from Equation \eqref{eq:defU}, the fourth  from $V=-W+\Psi$ and the fifth from $\Psi\indep W$ and $E(\Psi)=0$. The inequality holds since $(\sigma^2_C-\sigma^2_E)/(W^2+\sigma^2_C-\sigma^2_E) < 1$ with probability 1, as $\sigma^2_C>\sigma^2_E$.

\subsubsection*{Proof of Point 3}

Let $(V,W_h)$ be a normal vector with means $(0,h)$, variances $(\sigma^2_C,\sigma^2_C-\sigma^2_E)$, and covariance $-\left(\sigma^2_C-\sigma^2_E\right)$. Let $U_h=V+W_h(\sigma^2_C-\sigma^2_E)/(W_h^2+\sigma^2_C-\sigma^2_E)$. We have
\begin{align*}
&\sqrt{n}\left(\widehat{\beta}_{MSE}-\beta_0\right)\\
=&\sqrt{n}\left(\widehat{\beta}_{C}-\beta_0\right)\\
+&\frac{n\widehat{V}(\widehat{\beta}_C)-n\widehat{\cov}(\widehat{\beta}_C,\widehat{ \beta}_E)}{\left(\sqrt{n}\left(\widehat{\beta}_E-\beta_E-\left(\widehat{\beta}_C-\beta_0\right)\right)+h\right)^2+n\widehat{V}(\widehat{\beta}_E-\widehat{\beta}_C)}\left(\sqrt{n}\left(\widehat{\beta}_E-\beta_E-\left(\widehat{\beta}_C-\beta_0\right)\right)+h\right)\\
\convD&~U_h.
\end{align*}
The first equality follows from Equations \eqref{eq:diff_MSE_C} and  \eqref{eq:popt_closedform} and from $\beta_E=\beta_0+h/\sqrt{n}$. The convergence in distribution arrow follows from Assumption \ref{hyp:setup}, the Slutsky lemma, and the continuous mapping theorem.

\medskip
Let $\Psi_h=V+W_h$. The vector $(\Psi_h,W_h)$ is normally distributed, $\cov(\Psi_h,W_h)=0$ so $\Psi_h\indep W_h$. 
Let $g=h/\sqrt{\sigma^2_C-\sigma^2_E}$, and let $N_{g}=\frac{W_h}{\sqrt{\sigma^2_C-\sigma^2_E}}$, so that $N_{g}\sim\mathcal{N}(g,1)$. We have
\begin{align}\label{eq:comparison_MSE_local_alt}
&E(U_h^2)-\sigma^2_C \nonumber\\
=&E((U_h-V)(U_h+V))\nonumber\\
=&E\left(W_h\frac{\sigma^2_C-\sigma^2_E}{W_h^2+\sigma^2_C-\sigma^2_E}\left(2V+W_h\frac{\sigma^2_C-\sigma^2_E}{W_h^2+\sigma^2_C-\sigma^2_E}\right)\right)\nonumber\\
=&E\left(W_h\frac{\sigma^2_C-\sigma^2_E}{W_h^2+\sigma^2_C-\sigma^2_E}\left(-2W_h+2\Psi_h+W_h\frac{\sigma^2_C-\sigma^2_E}{W_h^2+\sigma^2_C-\sigma^2_E}\right)\right)\nonumber\\
=&-2\left(E\left(W_h\frac{\sigma^2_C-\sigma^2_E}{W_h^2+\sigma^2_C-\sigma^2_E}W_h\right)-E\left(W_h\frac{\sigma^2_C-\sigma^2_E}{W_h^2+\sigma^2_C-\sigma^2_E}\right)E\left(\Psi_h\right)\right)\nonumber\\
+&E\left(\left(W_h\frac{\sigma^2_C-\sigma^2_E}{W_h^2+\sigma^2_C-\sigma^2_E}\right)^2\right)\nonumber\\
=&-2\left(E\left(W_h\frac{\sigma^2_C-\sigma^2_E}{W_h^2+\sigma^2_C-\sigma^2_E}W_h\right)-E\left(W_h\frac{\sigma^2_C-\sigma^2_E}{W_h^2+\sigma^2_C-\sigma^2_E}\right)E\left(W_h\right)\right)\nonumber\\
+&E\left(\left(W_h\frac{\sigma^2_C-\sigma^2_E}{W_h^2+\sigma^2_C-\sigma^2_E}\right)^2\right)\nonumber\\
=&-2\cov\left(W_h\frac{\sigma^2_C-\sigma^2_E}{W_h^2+\sigma^2_C-\sigma^2_E},W_h\right)+E\left(\left(W_h\frac{\sigma^2_C-\sigma^2_E}{W_h^2+\sigma^2_C-\sigma^2_E}\right)^2\right)\nonumber\\
=&\left(\sigma^2_C-\sigma^2_E\right)\left\{E\left[\left(\frac{N_{g}}{N_{g}^2+1}\right)^2\right]-2\cov\left(\frac{N_{g}}{N_{g}^2+1},N_{g}\right)\right\}.
\end{align}
Let
\begin{equation}\label{eq:def_Delta(g)}
\Delta(g)=E\left[\left(\frac{N_{g}}{N_{g}^2+1}\right)^2\right]-2\cov\left(\frac{N_{g}}{N_{g}^2+1},N_{g}\right).
\end{equation}
Since $N_{-g}\sim -N_{g}$, we have
$$\Delta(-g)=E\left[\left(\frac{-N_{g}}{(-N_{g})^2+1}\right)^2\right]-2\cov\left(\frac{-N_{g}}{(-N_{g})^2+1},-N_{g}\right)=\Delta(g).$$
Moreover, 
we obtain through numerical simulations that
$\min_{g\in\R^+}\;\Delta(g)\simeq -0.53$ and $\max_{g\in\R^+}\;\Delta(g)\simeq 0.25$. The inequality
$$\max_{h\in \R}\left[(E(U_h))^2+V(U_h)-\sigma^2_C\right]<\max_{h\in \R}\left[\sigma^2_C-\left((E(U_h))^2+V(U_h)\right)\right]$$
follows from these last results, Equations \eqref{eq:comparison_MSE_local_alt} and $\sigma^2_C>\sigma^2_E$.


\subsection{Proof of Proposition \ref{prop:reduced_bias}} 
\label{sub:proof_of_proposition_ref_prop_reduced_bias}

\subsubsection*{Proof of Point 1}

First, we have
$$\widehat{p}_{\lambda}=\frac{\widehat{V}(\widehat{\beta}_C)-\widehat{\cov}(\widehat{\beta}_C,\widehat{\beta}_E)}{\max\left[0,\left(\widehat{\beta}_E-\widehat{\beta}_C\right)^2-\lambda\Vh (\widehat{\beta}_E-\widehat{\beta}_C) \right]+\widehat{V}(\widehat{\beta}_E-\widehat{\beta}_C)}.$$
Then, if $\beta_0\neq\beta_E$, Assumption \ref{hyp:setup} and the continuous mapping theorem yields
$$n\widehat{p}_{\lambda}\convP \frac{\sigma^2_C - \sigma^2_E}{(\beta_E-\beta_C)^2}.$$
The result follows as in the previous proof.

\subsubsection*{Proof of Point 2}

Instead of considering $U_0$ defined by \eqref{eq:defU}, we now consider $U_{0,\lambda}$ defined by
$$U_{0,\lambda} = V + W\frac{\sigma^2_C- \sigma^2_E}{\max\left[0,W^2- \lambda(\sigma^2_C- \sigma^2_E)\right] + \sigma^2_C- \sigma^2_E}.$$
Then, the proofs of the convergence in distribution, of $E(U_{0,\lambda})=0$ and $E\left(U_{0,\lambda}^2\right)<\sigma^2_C$ are identical to those above. We  now show that $\lambda\mapsto V(U_{0,\lambda})$ is decreasing. Let  $Z=W/(\sigma^2_C-\sigma^2_E)^{1/2}$ Using the same definition of $\Psi$ as above, we have
$$U_{0,\lambda} =\Psi +  (\sigma^2_C-\sigma^2_E)^{1/2} Z\left(\frac{1}{\max(0,Z^2-\lambda)+1} - 1\right),$$
with $\Psi\indep Z$. Then, for any two $\lambda>\lambda'$,
\begin{align*}
	 V(U_{0,\lambda})-V(U_{0,\lambda'}) = & E\left[(U_{0,\lambda} +U_{0,\lambda'})(U_{0,\lambda}- U_{0,\lambda'})\right] \\
	 = & (\sigma^2_C-\sigma^2_E) E\left[Z^2\left(\frac{1}{\max(0,Z^2-\lambda)+1} + \frac{1}{\max(0,Z^2-\lambda')+1} - 2 \right) \right.\\
	 & \left. \hspace{2.2cm}\times \left(\frac{1}{\max(0,Z^2-\lambda)+1} - \frac{1}{\max(0,Z^2-\lambda')+1} \right)\right].
\end{align*}
Moreover,
$$\frac{1}{\max(0,Z^2-\lambda)+1} + \frac{1}{\max(0,Z^2-\lambda')+1} - 2 \leq 0,$$
and the inequality is strict for $|Z|>\sqrt{\lambda'}$. Also,
$$\frac{1}{\max(0,Z^2-\lambda)+1} - \frac{1}{\max(0,Z^2-\lambda')+1} \geq 0,$$
and again, the inequality is strict for $|Z|>\sqrt{\lambda'}$. Hence, $V(U_{0,\lambda})<V(U_{0,\lambda'})$. Finally, as $\lambda\to+\infty$, $U_{0,\lambda}$ converges almost surely to $\Psi$. Moreover, $|U_{0,\lambda}|\leq |\Psi| +  (\sigma^2_C-\sigma^2_E)^{1/2} |Z|$, so by the dominated convergence theorem, $V(U_{0,\lambda})\to V(\Psi)=\sigma^2_E$.

\subsubsection*{Proof of Point 3}

With the same notation as above, let
$$U_{h,\lambda} = \Psi_h + W_h \left(\frac{\sigma^2_C - \sigma^2_E}{\max(0,W_h^2 - \lambda(\sigma^2_C - \sigma^2_E)) + \sigma^2_C - \sigma^2_E} - 1\right).$$
Using a reasoning similar to that in the proof of Point 3 of Theorem \ref{thm_main},
\begin{align*}
&\sqrt{n}\left(\widehat{\beta}_{MSE,\lambda}-\beta_0\right)\convD~U_{h,\lambda}.
\end{align*}
By the same reasoning as in the proof of Point 3 of Theorem \ref{thm_main},
\begin{align*}
E(U_{h,\lambda}^2)- E(U_{h,0}^2) = (\sigma^2_C-\sigma^2_E)
 & \left\{ E\left[\left(\frac{N_g}{\max(0,N_{g}^2-\lambda)+1}\right)^2  - \left(\frac{N_g}{N_{g}^2+1}\right)^2\right]\right.  \\
 & \; \left. -2\cov\left[N_{g}\left(\frac{1}{\max(0,N_{g}^2-\lambda)+1} - \frac{1}{N_g^2+1}\right),N_{g}\right]\right\}.	
\end{align*}
Then, simulations show that \eqref{eq:ineq_minimax} holds.


\subsection{Proof of Theorem \ref{thm_main_nonparam}} 
\label{sub:proof_of_theorem_ref}

\subsubsection*{Proof of Point 1}

We have
\begin{equation}\label{eq:diff_MSE_C2}
\widehat{\beta}_{MSE}-\widehat{\beta}_C=\widehat{p}\left(\widehat{\beta}_E-\widehat{\beta}_C\right).
\end{equation}
If $\beta_0\ne\beta_E$, it follows from Assumption \ref{hyp:setup_nonparam} and the continuous mapping theorem that
$$r_n^2\left(\widehat{\beta}_{MSE}-\widehat{\beta}_C\right)\convP\frac{\sigma^2_C}{\beta_E-\beta_0}.$$
Therefore,
$$r_n\left(\widehat{\beta}_{MSE}-\widehat{\beta}_C\right)=o_P(1).$$

\subsubsection*{Proof of Point 2}

Let $V$ be a normal variable with mean $\mu$ and variance $\sigma^2_C$.
Let
\begin{equation}\label{eq:defU_nonparam}
U_0=V\left(1-\frac{\sigma^2_C}{V^2+\sigma^2_C}\right).
\end{equation}
If $\beta_0=\beta_E$,
\begin{align*}
&r_n\left(\widehat{\beta}_{MSE}-\beta_0\right)\\
=&r_n\left(\widehat{\beta}_{C}-\beta_0\right)\\
+&\frac{r_n^2\widehat{V}(\widehat{\beta}_C)-r_n^2\widehat{\cov}(\widehat{\beta}_C,\widehat{ \beta}_E)}{\left(r_n\left(\widehat{\beta}_E-\beta_E-\left(\widehat{\beta}_C-\beta_0\right)\right)\right)^2+r_n^2\widehat{V}(\widehat{\beta}_E-\widehat{\beta}_C)}r_n\left(\widehat{\beta}_E-\beta_E-\left(\widehat{\beta}_C-\beta_0\right)\right)\\
\convD&~U_0.
\end{align*}
The first equality follows from Equations \eqref{eq:diff_MSE_C2} and  \eqref{eq:popt_closedform} and from $\beta_E=\beta_0$. The convergence in distribution arrow follows from Assumption \ref{hyp:setup_nonparam}, the Slutsky lemma, and the continuous mapping theorem. Further,
\begin{align*}
E\left(U_0^2\right)-E\left(V^2\right)=&E((U_0-V)(U_0+V))\\
=&E\left(-\frac{\sigma^2_C}{V^2+\sigma^2_C}V\left(2V-\frac{\sigma^2_C}{V^2+\sigma^2_C}V\right)\right)\\
=&E\left(-\frac{\sigma^2_C}{V^2+\sigma^2_C}V^2\left(2-\frac{\sigma^2_C}{V^2+\sigma^2_C}\right)\right)\\
<&0.
\end{align*}
The third equality follows from Equation \eqref{eq:defU_nonparam}. The inequality follows from the fact
$1>\frac{\sigma^2_C}{V^2+\sigma^2_C}$
with probability 1.

\subsubsection*{Proof of Point 3}

Let
\begin{equation}\label{eq:defUh_nonparam}
U_h=V + (h-V)\frac{\sigma^2_C}{(h-V)^2+\sigma^2_C}.
\end{equation}
\begin{align*}
&r_n\left(\widehat{\beta}_{MSE}-\beta_0\right)\\
=&r_n\left(\widehat{\beta}_{C}-\beta_0\right)\\
+&\frac{r_n^2\widehat{V}(\widehat{\beta}_C)-r_n^2\widehat{\cov}(\widehat{\beta}_C,\widehat{ \beta}_E)}{\left(r_n\left(\widehat{\beta}_E-\beta_E-\left(\widehat{\beta}_C-\beta_0\right)\right)+h\right)^2+r_n^2\widehat{V}(\widehat{\beta}_E-\widehat{\beta}_C)}\left(r_n\left(\widehat{\beta}_E-\beta_E-\left(\widehat{\beta}_C-\beta_0\right)\right)+h\right)\\
\convD&~U_h.
\end{align*}
The first equality follows from Equations \eqref{eq:diff_MSE_C2} and  \eqref{eq:popt_closedform} and from $\beta_E=\beta_0+h/r_n$. The convergence in distribution arrow follows from Assumption \ref{hyp:setup_nonparam}, the Slutsky lemma, and the continuous mapping theorem.

\medskip
Let $g=\frac{h-\mu}{\sigma_C}$, $N_g=\frac{h-V}{\sigma_C}$, and $\mu_{sd}=\frac{\mu}{\sigma_C}$. Then $N_g\sim \mathcal{N}(g,1)$ and we have:
\begin{align}\label{eq:comparison_MSE_local_alt_non_param}
E\left(U_h^2\right)-E\left(V^2\right)=&E((U_h-V)(U_h+V))\nonumber\\
=&E\left((h-V)\frac{\sigma^2_C}{(h-V)^2+\sigma^2_C}\left(2V+(h-V)\frac{\sigma^2_C}{(h-V)^2+\sigma^2_C}\right)\right)\nonumber\\
=&\sigma^2_CE\left(\frac{N_g}{N_g^2+1}\left(2\left(g+\mu_{sd}-N_g\right)+\frac{N_g}{N_g^2+1}\right)\right)\nonumber\\
=&\sigma^2_C\left(\Delta(g)+2\mu_{sd}E\left(\frac{N_g}{N_g^2+1}\right)\right),
\end{align}
with $\Delta(g)$ defined as in Equation \eqref{eq:def_Delta(g)}.
Let $\Lambda(g,\mu_{sd})=\Delta(g)+2\mu_{sd}E\left(\frac{N_g}{N_g^2+1}\right)$. The function $x\mapsto x/(x^2+1)$ and the density of a $\mathcal{N}(0,1)$ satisfy the conditions of Lemma \ref{lem:expect_symm}. Thus, $E\left(\frac{N_g}{N_g^2+1}\right)\geq 0$ if $g\geq 0$, and $E\left(\frac{N_g}{N_g^2+1}\right)<0$ otherwise. Then, if $\mu$ and $g$ are not of the same sign, $\Lambda(g,\mu_{sd})\leq \Delta(g)$, so it follows from Point 3 of Theorem \ref{thm_main} that for every $\mu_{sd}$, the minimum of $\Lambda(g,\mu_{sd})$ with respect to $g$ is greater in absolute value than its maximum. We can then restrict attention to cases where $\mu$ and $g$ are of the same sign. As $\Lambda(-g,-\mu_{sd})=\Lambda(g,\mu_{sd})$, we can further restrict attention to cases where $g\geq 0$ and $\mu_{sd}\geq 0$. We use Monte-Carlo simulations with $10^7$ Halton draws to approximate $\Lambda(g,\mu_{sd})$ for every $g\in\{0,0.1,0.2,...,10\}$ and  $\mu_{sd}\in[0,0.41]$, 
the absolute value of the minimum of $\Lambda(g,\mu_{sd})$ with respect to $g$ is larger than the absolute value of its maximum. This proves the result, together with Equation \eqref{eq:comparison_MSE_local_alt_non_param}.





\newpage
\bibliography{biblio}

\end{document}